\documentclass[letterpaper, 10 pt, conference]{ieeeconf}  % Comment this line out if you need a4paper

\IEEEoverridecommandlockouts                              % This command is only needed if 
\usepackage{tikz}
\usepackage{graphics} % for pdf, bitmapped graphics files
\usepackage{epsfig} % for postscript graphics files
\usepackage{amsmath} % assumes amsmath package installed
\usepackage{amssymb}  % assumes amsmath package installed
\usepackage{url,algorithm}
\usepackage{multirow}
\usepackage{color}
\usepackage{adjustbox}
\usepackage{booktabs}
\usepackage{array}

\let\labelindent\relax
\usepackage{enumitem}
\renewcommand{\theenumi}{\arabic{enumi}}

\renewcommand{\theenumii}{\arabic{enumii}}

\renewcommand{\theenumiii}{\arabic{enumiii}}

\renewcommand{\theenumiv}{\arabic{enumiv}}

{\normalfont}{\normalfont}

\newcommand{\Nu}{N_u}

\newcommand{\Np}{N_p}

\newcommand{\parvec}{x}

\newcommand{\sizeparvec}{n_\parvec}
\newcommand{\pref}{\pi}
\newcommand{\rr}{{\mathbb R}}
\newcommand{\domain}{\mathcal{D}}
\newcommand{\domainHARD}{\Omega_G}
\newcommand{\domainSOFT}{\Omega_S}
\newcommand{\ba}[1]{\begin{array}{#1}}
	\newcommand{\ea}{\end{array}}
\newcommand{\st}{\mathop{\rm s.t.}\nolimits}

\newcommand{\surr}[1]{\hat{#1}}

\def\BibTeX{{\rm B\kern-.05em{\sc i\kern-.025em b}\kern-.08em
    T\kern-.1667em\lower.7ex\hbox{E}\kern-.125emX}}
\begin{document}
\title{C-GLISp: Preference-Based Global Optimization under Unknown Constraints
with Applications to Controller Calibration}
\author{Mengjia Zhu, Dario Piga, and Alberto Bemporad
\thanks{This paper was partially supported by the Italian Ministry of University and Research under the PRIN'17 project ``Data-driven learning of constrained control systems'' , contract no. 2017J89ARP.}
\thanks{M. Zhu and A. Bemporad are with IMT School for Advanced Studies Lucca, 55100 Lucca, Italy (e-mail: mengjia.zhu@imtlucca.it; alberto.bemporad@imtlucca.it).}
\thanks{D. Piga is with IDSIA Dalle Molle Institute for Artificial Intelligence, SUPSI-USI, 6962 Lugano, Switzerland (e-mail: dario.piga@supsi.ch).}}

\maketitle

\begin{abstract}
Preference-based global optimization algorithms minimize
an unknown objective function only based on
whether the function is better, worse, or similar for given pairs of candidate optimization vectors. Such optimization problems arise in many real-life examples, such as finding the optimal calibration of the parameters of a control law. The calibrator can judge whether a particular combination of parameters leads to a better, worse, or  similar closed-loop performance. 
Often, the search for the optimal parameters is also subject to unknown constraints. For example, the vector of calibration parameters must not lead to closed-loop instability.
This paper extends an active preference learning algorithm introduced recently by the authors to handle unknown constraints. 
The proposed method, called C-GLISp, looks for an optimizer of the problem
only based on \emph{preferences} expressed on pairs of candidate vectors, 
and on whether a given vector is reported \emph{feasible} and/or 
\emph{satisfactory}. C-GLISp learns a surrogate of the underlying objective function based on the expressed preferences, and a surrogate of the probability that a sample is feasible and/or satisfactory based on whether each of the tested vectors was judged as such.
The surrogate functions are used iteratively to propose a new candidate vector to test and judge. Numerical benchmarks and a semi-automated control calibration task
demonstrate the effectiveness of C-GLISp, showing that it can reach near-optimal solutions within a small number of iterations.

\end{abstract}

%\begin{IEEEkeywords}
%Active preference learning, Global Optimization with Unknown Constraints, Model predictive control
%\end{IEEEkeywords}

\section{Introduction}
\label{sec:introduction}
Active learning algorithms for black-box global optimization problems have been studied since the sixties under different names~\cite{Mat63,Kus64,SWMW89,JSW98}. These algorithms solve the problem by optimizing a surrogate of the objective function, which is estimated by exploring the space of the optimization variables.
In particular, nowadays \emph{Bayesian Optimization} (BO)~\cite{shahriari2015}
is widely used to solve problems in which the cost function can only be quantified 
after running an experiment, such as in experimental controller calibration~\cite{savaia2021experimental} and in automated machine learning~\cite{Luo16}.
The main idea of such methods is to fit a surrogate function to the available observations and iteratively suggest the next query point by optimizing an acquisition function. The latter trades off between exploiting the surrogate function to improve the objective and scouting unexplored areas 
of the search domain. An alternative surrogate method to BO based on estimating the underlying
objective by radial basis functions and inverse distance weighting for exploration,
called GLIS, was recently proposed in~\cite{Bem20}.

Successful applications of global optimization algorithms based on active learning  for the calibration of Model Predictive Control (MPC), PID, and state-feedback control laws were presented in~\cite{FoPiBe20,LuFoCo20,DBLP:journals/corr/abs-1906-12086,marco2016automatic}.
In this context, the tuning parameters of the controller are the optimization variables,
and a quantitative characterization of the resulting closed-loop performance after running a simulation or experiment is the objective to optimize. These algorithms were also used for model selection~\cite{piga2019performance,bansal2017goal},  for controller tuning in robotic  manipulation  and trajectory tracking~\cite{driess2017constrained,ROVEDA2021103711,roveda2020two}, in optimizing gait parameters in robotic  bipedal locomotion~\cite{calandra2014bayesian}, and for ``safe'' optimization of position controller parameters of quadrotors~\cite{berkenkamp2016safe}.  

A limitation of black-box optimizers like BO and GLIS is that they require quantifying an objective function
after running an experiment. However, many real-world controller calibration problems
involve multiple objectives to optimize, such as settling time, overshoots, actuation effort,
computational burden, and other performance-related metrics. The relative weights of such objectives can be hard 
to assign, and sometimes even impossible to quantify, as they are the result of a qualitative judgment. On the other hand, a skilled calibrator can often assess closed-loop performance of certain tuning combinations in terms of ``this test was better than the other one,''
\emph{i.e.}, by pairwise comparisons. Thus, when quantifying an objective function is difficult or impossible, one can instead use these expressed \emph{preferences} to learn an underlying surrogate function to  be optimized, which leads to the area of preference-based learning algorithms.

Preference-based Bayesian optimization (PBO) has been proposed in~\cite{chu2005extensions, brochu2007active, gonzalez2017preferential, abdolshah2019multi}. Preference-based \emph{reinforcement learning} (RL) has also drawn much attention in recent years~\cite{christiano2017deep}. The reader is referred to the survey paper~\cite{wirth2017survey} for a comprehensive review. Note that sample efficiency, which is related to the credit assignment task in RL, is a major challenge in many preference-based RL methods~\cite{wirth2017survey}. A more sample-efficient active preference learning method, called GLISp (an extension of GLIS), was proposed in~\cite{BemPig20}.
GLISp  learns a surrogate of the underlying preference relations by solving a Quadratic Programming (QP) problem, whose constraints reflect the expressed \emph{preference} on whether a specific candidate is better than the other. Then, the
algorithm iteratively proposes a new candidate for testing to the decision-maker for comparison. This experiment-driven preference-based approach was tested in~\cite{zhu2020pref} for semi-automated MPC calibration (automatic selection of  control parameters, manual  assessment of performance by comparisons), demonstrating its effectiveness in terms of the number of experiments needed to reach near-optimal closed-loop performance.

Real-life control design problems often involve constraints that are unknown beforehand, or for which it is not possible to find an explicit analytic expression. This is a challenge since safe exploration can be essential in many control applications, and infeasible experiments can be dangerous and costly. Several methods have been proposed in the literature to handle unknown constraints and encourage safe exploration. In~\cite{sui2018stagewise}, Sui \emph{et al.} presented a stagewise safe BO with Gaussian processes, which they later extended to allow multiple safety constraints independent of the objective function~\cite{berkenkamp2020bayesian}. A general formulation for constrained BO and a modified version of the expected improvement acquisition function was illustrated in~\cite{gelbart2014bayesian}, which  handles noisy constraint observations and considers cases in which the objective and constraint functions are decoupled. In~\cite{antonio2019sequential}, a \emph{sequential model-based optimization method} was proposed. The unknown feasible region boundaries are first reconstructed from data  through \emph{support vector machines}. Then, a global optimization step is performed via BO. Differently from the aforementioned methods, GLISp accounts for unknown constraints \emph{implicitly} in the preferences expressed by the decision maker by making the samples that are infeasible lose the pairwise comparisons against the feasible ones. 

This paper extends GLISp to handle   
unknown constraints in the active learning phase \emph{explicitly}, therefore encouraging safe exploration. 
Besides expressing preferences, the decision-maker is asked to label an experiment
as feasible and if its outcome is overall satisfactory (yes/no).  Based on such labels, a surrogate of the probability of constraint feasibility and experiment's  satisfaction is learned via an \emph{Inverse Distance Weighting} (IDW) interpolant function~\cite{Bem20}. The surrogates are properly  integrated within  the acquisition function  to find the next point to test. 
We show the efficiency and effectiveness of the proposed method, called C-GLISp,  in three numerical benchmarks, and on an extension of the case study originally proposed in~\cite{zhu2020pref} on semi-automated MPC calibration for autonomous driving.  
MATLAB and Python implementations of C-GLISp are also provided and available at \url{http://cse.lab.imtlucca.it/~bemporad/glis}.

The rest of the paper is organized as follows. The problem of preference-based optimization with unknown constraints is formulated in Section~\ref{sec:problem}. The proposed active-learning algorithm and  details for its practical implementation are discussed in Section~\ref{sec:C-GLISp_formulation}.  Numerical benchmarks showing the properties and the effectiveness of the proposed method are reported  in Section~\ref{sec:benchmarks}, while the case study on   semi-automated MPC calibration for autonomous driving is presented in   Section~\ref{sec:case_study}. Conclusions and directions for future research are drawn in Section~\ref{sec:conclusion_futureExtension}.

\section{Problem Formulation}\label{sec:problem}
Let $\domain \subseteq \rr^{\sizeparvec}$ be the space of decision vectors $\parvec$. We are interested in minimizing an (unknown) objective function $f: \domain \rightarrow \mathbb{R}$ subject to the constraint that  $x$ belongs to an (unknown) feasibility set $\domainHARD \subseteq \domain$.

We  assume that we cannot represent the set $\domainHARD$, but rather that,
given a vector $\parvec \in \domain$, a decision-maker
can assess the value of the \emph{feasibility function} $G: \domain\to\{0,1\}$ defined as 
\begin{equation}
G(\parvec)=\left\{\ba{ll}
0 &\mbox{if $\parvec\notin\domainHARD$ }\\
1 &\mbox{if $\parvec\in\domainHARD$}.
\ea\right.
\label{eq:fes_fun}
\end{equation}
In other words, a value $\parvec$ outside  $\domainHARD$  will be considered as ``unacceptable'' ($G(x)=0$) by the decision-maker. For example, an unacceptable $\parvec$ can be a set of  controller parameters   leading to an unstable closed-loop behavior or to a control law that is too expensive
to compute in real-time.  

Furthermore, we assume that the  objective  function $f$ cannot be directly quantified,
but rather that can be indirectly observed in two ways:
\begin{enumerate}
	\item For a given sample $x \in \domain$, the decision-maker is  requested to say whether or not  $x$  leads to certain ``satisfactory performance''. Formally, we can define a  \emph{satisfaction set} $\domainSOFT \subseteq \domain$ and a \emph{satisfaction function} $S: \domain\to\{0,1\}$ as
\begin{equation}
S(\parvec)=\left\{\ba{ll}
0 &\mbox{if $\parvec \notin \domainSOFT$}\\
1 &\mbox{if $\parvec \in \domainSOFT$},
\ea\right.
\label{eq:softcost_fun}
\end{equation}
where the set $\domainSOFT$ contains all the vectors $\parvec$ leading to a performance that
the decision-maker judges satisfactory. An analytic  expression of $\domainSOFT$ is therefore
not available, only the value $S(x)$ is provided by the decision-maker
for any given $x \in \domain$. Note that $\domainSOFT$ may not be a subset of $\domainHARD$, for example when the preference-based optimization process is carried out in simulation: a sample may lead to satisfactory performance but would not be implementable due to hardware limitations.  
On the other hand, in cases of assessments based on physical experiments, $\domainSOFT$ is necessarily a subset of $\domainHARD$, as no performance would be available for evaluation when the parameters are infeasible. %satisfactory performance would not be possible to achieve if the parameters are infeasible.
%\item   For any  input pair  $\parvec_1, \parvec_2 \in \domainHARD$,
\item   For any pair  $\parvec_1,\parvec_2 \in \domain$, the decision-maker is  requested to provide  the output of the   \emph{preference function}:  $\pref:\domain\times\domain\to\{-1,0,1\}$ 
\begin{equation}
\pref(\parvec_1,\parvec_2)=\left\{\ba{ll}
-1 &\mbox{if $\parvec_1$ ``better'' than $\parvec_2$}\\
0 &\mbox{if $\parvec_1$ ``as good as'' $\parvec_2$}\\
1 &\mbox{if $\parvec_2$ ``better'' than $\parvec_1$}.
\ea\right.
\label{eq:pref_fun}
\end{equation}
where the preference in~\eqref{eq:pref_fun} is implicitly defined according to the 
underlying hidden function $f$ to be minimized, namely
\begin{equation}
\pref(\parvec_1,\parvec_2)=\left\{\ba{ll}
-1 &\mbox{if $  f(\parvec_1)<  f(\parvec_2)$}\\
0 &\mbox{if $  f(\parvec_1)=  f(\parvec_2)$}\\
1 &\mbox{if $  f(\parvec_1)>  f(\parvec_2)$}. 
\ea\right.
\label{eq:pref_fun-f}
\end{equation}
\end{enumerate}

The rationale behind the above problem formulation is that often one encounters practical decision problems in which a function $f$ is impossible to quantify, but anyway it is possible for a human operator to express a qualitative evaluation (\emph{e.g.}, ``good'' or ``bad'')  and 
a preference between the outcome of two experiments. 
 
Formally, we want to find the optimal solution $\parvec^{\star} \in \domainSOFT \cap  \domainHARD$ such that $\parvec^{\star}$
is ``better'' (or ``no worse'') than any other $\parvec$ according to the preference function $\pref$:
\begin{equation}
\mbox{find}\ \parvec^\star\ \mbox{such that}\ \pref(\parvec^\star,\parvec)\leq 0,\ \forall \parvec \in \domainSOFT \cap  \domainHARD.
\label{eq:glob-opt-pref}
\end{equation}
We propose to solve problem~\eqref{eq:glob-opt-pref} iteratively as follows: ($i$) suggest a sequence of  decision vectors $\parvec_1,\ldots,\parvec_N \in \domain$ to test, ($ii$) ask to evaluate the feasibility function $G(x_i)$ and the satisfaction function $S(x_i)$ for $i=1,\ldots,N$, and ($iii$) ask to evaluate the preference function  $\pi(x_i,x_j)$ for $M$ given pairs $(i,j)$, $i,j=1,\ldots,N$, $i\neq j$, where $M$ is the number of expressed preferences, $1 \leq M \leq \binom{N}{2}$. The goal is to propose candidate vectors $\parvec_N$ approaching the optimal solution $\parvec^\star$ as $N$ grows.

\section{Proposed Method} \label{sec:C-GLISp_formulation}
The proposed preference-based optimization method to solve problem~\eqref{eq:glob-opt-pref} is based on an extension of the GLISp 
algorithm originally introduced in~\cite{BemPig20}. We refer to the new algorithm
as C-GLISp, whose aim is to handle unknown constraints expressed in terms of an approximation
of the feasibility function $G$  in~\eqref{eq:fes_fun} and of the satisfaction function $S$ in~\eqref{eq:softcost_fun}. 

Similarly to GLISp, C-GLISp involves two main phases: an initial random sampling  and an active learning phase. In both phases, C-GLISp trains and updates three surrogate functions approximating, respectively, the feasibility function $G$,  the satisfaction function $S$, and the underlying function $f$. During the active learning phase, the next point for evaluation is selected by optimizing an acquisition function which trades off \emph{exploitation} (optimization only based on the surrogates describing the observed preferences and constraints) and  \emph{exploration}  (searching unexplored areas of the domain $\domain$). The goal of C-GLISp is to approach an optimal solution 
% [AB: the global optimizer may not be unique]
$\parvec^\star$ as in~\eqref{eq:glob-opt-pref} within a small number $N$ of experiments. %~\cite{BemPig20,zhu2020pref}. 

\subsection{Learning Unknown Constraint Functions}
We discuss how to train surrogates of the functions $G$ %~\eqref{eq:fes_fun} 
and $S$ %~\eqref{eq:softcost_fun} 
that approximate, respectively, the feasibility constraint  $\parvec \in \domainHARD$ and the satisfaction constraint $\parvec \in \domainSOFT$. 
The idea is to ask the decision-maker to assess, once an experiment is performed, whether the constraints  $\parvec \in \domainHARD$ and $\parvec \in \domainSOFT$ are satisfied or not,
and train surrogate functions of $G$ and $S$ based on the outcome of $N\geq 2$ of such queries.
These queries  are performed on a set of samples $\{\parvec_1,\ldots,\parvec_N\}$ iteratively proposed by C-GLISp.

Compared to unconstrained preference-based optimization like GLISp,
in which an infeasible/unsatisfactory sample only indirectly reveals itself
as such by losing pairwise comparisons against feasible/satisfactory ones, C-GLISp exploits the information 
on whether $\parvec \in \domainHARD$ and/or $\parvec \in \domainSOFT$ to facilitate the optimization process, in particular, to avoid exploring the infeasible and/or unsatisfactory region and therefore reduce the number of samples $\parvec_i\not \in \domainHARD$ and/or $\parvec_i \not \in \domainSOFT$.

The surrogate functions for $G$ and  $S$ are constructed as follows. The decision-maker observes  the outcome of the performed experiments, and he/she provides a \emph{feasibility vector} $G_F =[G_1\ \ldots\ G_N]'\in\{0,1\}^{N}$ with 
\begin{equation}
G_i=G(x_i),
\label{eq:feasibility-vector}
\end{equation}
and a \emph{satisfaction vector} $S_F =[S_1\ \ldots\ S_N]'\in\{0,1\}^{N}$ with
\begin{equation}
S_i=S(x_i),
\label{eq:satisfactory-vector}
\end{equation}
by assessing whether each experiment is feasible and satisfactory. Then, surrogates $\surr{G}$ of $G$
and $\surr{S}$ of $S$ are constructed from the observations $G_F$ and $S_F$, respectively, as detailed below. 

A surrogate function $\surr{G}:\domain\to\rr$ predicting the probability of satisfying the feasibility constraint $\parvec \in \domainHARD$ is defined as 
\begin{equation}
\surr{G}(\parvec)=\sum_{i=1}^N\nu_i(\parvec)G_i,
\label{eq:idw_intep}
\end{equation}
where $\nu_i(x): \domain\to\rr $ for $i = 1\ldots,N$ is defined as
\begin{equation}
\nu_i(\parvec)=\left\{\ba{lll}
1 &\mbox{if $\parvec = \parvec_i$ }\\
0 &\mbox{if $\parvec = \parvec_j, j\neq i$}\\
\frac{w_i(\parvec)}{\sum_{i=1}^Nw_i(\parvec)} &\mbox{otherwise}. 
\ea\right.
\label{eq:idw_coef_nu}
\end{equation}

Here $w_i:\domain \setminus \{\parvec_i\} \rightarrow \mathbb{R}$ is the following IDW function~\cite{joseph2011regression} 
\begin{equation}
    w_i(x)=\frac{e^{-d^2(x,x_i)}}{d^2(x,x_i)},
\label{eq:w-IDW}
\end{equation}
where  $d:\domain\times\domain\to\rr$  denotes  the squared
Euclidean distance
\begin{equation}
d(\parvec,\parvec_i)=\|\parvec-\parvec_i\|_2^2. 
\label{eq:distance}
\end{equation}

The surrogate function $\surr{S}:\domain\to\rr$ is defined similarly. The approach presented in this paper aims at solving problems where experiments are expensive to run, so that data efficiency is essential. IDW interpolation functions are selected
in this case because of their high accuracy. Other binary classification methods (\emph{e.g.} logistic regression or random forests) would be less suitable in this context since their accuracy with a small number of training data is limited. Support vector machines (SVMs)~\cite{SS04} can be a potential substitute since they work well with small and medium-size training sets. However, our numerical tests have shown that IDW interpolation functions outperform SVM. In addition, the functions $\surr{G}$ and $\surr{S}$ generated by IDW interpolation are always between 0 and 1 by construction (see~\cite[Lemma 1-P2]{Bem20}), and can be interpreted as
probabilities of being feasible/satisfactory. 
%It is also worth noting that in case of multiple constraints, as they are unknown it may be impossible to distinguish between them. Hence, rather than modeling each constraint with a separate surrogate function, we model the probability of being feasible/satisfactory, making it easier to handle multiple feasibility/satisfactory constraints.

\subsection{Learning The Preference Function}
Radial basis functions (RBFs)~\cite{Gut01,MGTM07} are flexible and have been adopted to solve global optimization problems in~\cite{Bem20, costa2018rbfopt, Gut01, regis2005constrained,MGTM07} with success. Therefore, as in~\cite{BemPig20}, we parameterize the surrogate function $\hat f:\domain\to\rr$
as a linear combination of RBFs~\cite{Gut01,MGTM07}
\begin{equation}
\hat f(\parvec)=\sum_{k=1}^N\beta_k\phi(\epsilon d(\parvec,\parvec_i)),
\label{eq:rbf}
\end{equation}
where $\phi:\rr\to\rr$ is an RBF, 
$\epsilon>0$ is a scalar hyper-parameter defining the shape of the RBF, 
and $\beta=[\beta_1\ \ldots\ \beta_N]^{T}$ are the unknown coefficients to be determined. 
Some RBFs commonly used are $\phi(\epsilon d)=\frac{1}{1+(\epsilon d)^2}$
(\emph{inverse quadratic}), $\phi(\epsilon d)=e^{-(\epsilon d)^2}$ (\emph{Gaussian}), and $\phi(\epsilon d)=(\epsilon d)^2\log(\epsilon d)$ (\emph{thin plate spline}).

Besides the \emph{feasibility vector} $G_F$ and the satisfaction vector $S_F$, the \emph{preference vector} $B=[b_1\ \ldots\ b_{M}]^{T}\in\{-1,0,1\}^{M}$ is also assumed to be provided by the decision-maker, with
\begin{equation}
b_h=\pref(\parvec_{i(h)},\parvec_{j(h)}),
\label{eq:pref-vector}
\end{equation}
for $\parvec_i,\ \parvec_j\in\domain$ such that $\parvec_i\neq \parvec_j$, $\forall i\neq j$,$\ i,j=1,\ldots,N$. 
Here, $M$ is the number of expressed preferences, $1 \leq M \leq \binom{N}{2}$, 
$h\in\{1,\ldots,M\}$ is the index enumerating the preferences, and  $i(h), j(h)\in \{1,\ldots,N\}$, $i(h) \neq j(h)$. 

The preferences $b_h$ expressed by the decision-maker are used to shape the surrogate objective function $\hat f$ by imposing the following constraints:
\begin{equation}
\ba{ll}
\hat f(\parvec_{i(h)})\leq\hat f(\parvec_{j(h)})-\sigma& \mbox{if}\ b_h=-1\\
\hat f(\parvec_{i(h)})\geq\hat f(\parvec_{j(h)})+\sigma& \mbox{if}\ b_h=1\\
|\hat f(\parvec_{i(h)})-\hat f(\parvec_{j(h)})|\leq \sigma& \mbox{if}\ b_h=0
\ea
\label{eq:RBF-pref}
\end{equation}
for $h=1,\ldots,M$, where $\sigma>0$ is a given scalar
that avoids the trivial solution $\hat f(x)\equiv 0$. 

Similarly to SVMs~\cite{SS04}, the vector $\beta$ of coefficients describing the surrogate $\hat f$  is obtained by solving the following convex QP problem 
%support vector machine (SVM)
\begin{equation}
\ba{rll}
\min_{\beta,\varepsilon} &\displaystyle{\sum_{h=1}^{M} c_h\varepsilon_h+\frac{\lambda}{2}\sum_{k=1}^N\beta_{k}^2 }\\
\st % & \sum_{j=1}^N \phi(\epsilon d(x_{t1},x_{j}))\beta_j=0\\
% &\sum_{j=1}^N \phi(\epsilon d(x_{s2},x_{j}))\beta_j=1\\
& \displaystyle{\sum_{k=1}^N\beta_k(\phi(\epsilon d(\parvec_{i(h)},\parvec_{k}))-\phi(\epsilon d(\parvec_{j(h)},\parvec_{k})))}
\\&\quad\quad\leq -\sigma+\varepsilon_h,\hfill \forall h:\ b_h=-1\\
&\displaystyle{\sum_{k=1}^N\beta_k(\phi(\epsilon d(\parvec_{i(h)},\parvec_{k}))-\phi(\epsilon d(\parvec_{j(h)},\parvec_{k})))}\\
&\quad\quad\geq \sigma-\varepsilon_h, \hfill \forall h:\ b_h=1\\
&\displaystyle{\left| \sum_{k=1}^N\beta_k(\phi(\epsilon d(\parvec_{i(h)},\parvec_{k}))-\phi(\epsilon d(\parvec_{j(h)},\parvec_{k})))\right|}\\&
\quad\quad\leq\sigma+\varepsilon_h,\hfill \forall h:\ b_h=0\\  & h=1,\ldots,M
\ea
\label{eq:QP-pref}
\end{equation}
that captures the preference relationships in~\eqref{eq:RBF-pref} and the parametrization of $\hat f$ in~\eqref{eq:rbf}. In~\eqref{eq:QP-pref}, 
$c_h$ are positive weights and $\varepsilon_h$ are positive slack variables  used to relax the constraints imposed by~\eqref{eq:RBF-pref}. The violation of the imposed constraints could be caused by an inappropriate selection of the RBF, leading to poor flexibility in the parametric description of the surrogate function $\hat f$, and by inconsistent assessments provided by the decision-maker. The scalar $\lambda>0$ in the cost function~\eqref{eq:QP-pref} is a regularization parameter. With $\lambda>0$, problem~\eqref{eq:QP-pref} is a QP problem that admits a unique solution.

In the constrained preference-based optimization algorithm detailed in the following section, we will execute $K$-fold cross-validation periodically (\emph{i.e.}, when $i$ is in the predefined self-calibration index set  $\mathcal{I}_{\rm{sc}} \subseteq \{1,\ldots,N_{\rm max}-1\}$) to automatically tune the hyper-parameter $\epsilon$ defining the shape of the RBF in~\eqref{eq:rbf} during the active learning phase, as recommended in~\cite{BemPig20}. 

\subsection{Acquisition Function}\label{sec:acquisition}
Minimizing $\hat f$ greedily to generate the next sample $\parvec_{N+1}$ may lead the solver to converge to a point that is not the global optimum of~\eqref{eq:glob-opt-pref}. Hence, when selecting the next point $\parvec_{N+1}$, besides \emph{exploiting} the surrogate $\hat f$, some \emph{exploration} should be considered to search regions with limited/no samples to reduce the uncertainty associated with $\hat f$. Also, the feasibility and satisfactory regions are unknown and are only \emph{implicitly} included in the surrogate function $\hat f$. Therefore, we also include terms to \emph{explicitly}  avoid the \emph{exploration} in the regions with low probabilities of being feasible and satisfactory by penalizing the (estimated) infeasibility $x\not\in\domainHARD$ and unsatisfactory performance $\parvec\not\in \domainSOFT$.

The exploration term used in GLISp is the following IDW function 
$z:\domain\to\rr$
\begin{equation}
z(\parvec)=\left\{\ba{ll}
0 & \mbox{if}\ \parvec\in\{\parvec_1,\ldots,\parvec_N\}\\
%\frac{2}{\pi}
\tan^{-1}\left(\frac{1}{ { \sum_{i=1}^N r_i(x) }}\right)&
\mbox{otherwise,}\ea\right.
\label{eq:IDW-distance}
\end{equation}
where $r_i(x)=\frac{1}{d^2(\parvec,\parvec_i)}$. 
Note that $z(\parvec)=0$ for all the decision variables $\parvec$ already sampled and tested, and $z(\parvec)>0$ in $\domain\setminus \{\parvec_1,\ldots,\parvec_N\}$. The \emph{arc tangent} function is used to prevent the new sampled point from getting excessively far away from the existing ones.

Unlike GLISp, here we modify~\eqref{eq:IDW-distance} into
\begin{equation}
\begin{split}
z_N(\parvec) = &\left(1 - \frac{N}{N_{max}}\right) \tan^{-1}\left(\frac{\sum_{i=1}^Nr_i(\parvec^*_N)}{\sum_{i=1}^Nr_i(\parvec)}\right) \\
 &+\frac{N}{N_{max}}\tan^{-1}\left(\frac{1}{\sum_{i=1}^Nr_i(\parvec)}\right)
\end{split}
\label{eq:z_comb}
\end{equation}
for $x \notin\{\parvec_1,\ldots,\parvec_N\}$ and $z_N(x) = 0$ otherwise. In~\eqref{eq:z_comb}, $x_N^*$ is the best decision variable found up to iteration $N$.
In~\eqref{eq:z_comb}, $N_{max}$ is the maximum allowed number of experiments.  The rationale behind~\eqref{eq:z_comb} is that it encourages the exploration of regions of $\domain$ further away from the current best solution in the early iterations and reduce its effects as the number $N$ of experiments increases.   The exploration function in~\eqref{eq:z_comb} was empirically observed to better escape from local minima than that in~\eqref{eq:IDW-distance}.

The \emph{acquisition function} $a:\domain\to\rr$ is defined as
\begin{equation}
	\begin{split}
a(\parvec)= &\frac{\hat f(\parvec)}{\Delta\hat F}-\delta_E z_N(\parvec) \\
&+ \delta_G(1- \surr{G}(\parvec)) +  \delta_S(1- \surr{S}(\parvec)),
\end{split}
\label{eq:acquisition_modified}
\end{equation}
where $\delta_E \geq 0$ is the exploration parameter, and $\delta_G,\delta_S \geq 0$ weight the probability of a sample $\parvec$ to be infeasible and/or unsatisfactory, respectively. 
The term $\Delta \hat F=\max_i\{\hat f(\parvec_i)\}-\min_i\{\hat f(\parvec_i)\}$ 
is the range of the surrogate function $\hat f$ on the samples in $\{\parvec_1,\ldots,\parvec_N\}$.
It is used as a scaling factor in~\eqref{eq:acquisition_modified} to make each term in~\eqref{eq:acquisition_modified} comparable, which eases the selection of the hyper-parameters $\delta_E$, $\delta_G$, and $\delta_S$. % When there is at least one comparison $b_h=\pref(\parvec_{i(h)},\parvec_{j(h)}) \neq 0$, $\Delta \hat F \geq \sigma$ .

The exploration parameter $\delta_E$ encourages sampling unexplored regions of the domain $\domain$. Setting $\delta_E=0$ makes C-GLISp rely heavily on the accuracy of the surrogate functions $\hat f$, $\hat G$, and $\hat S$, which may easily lead to missing the global  optimum. On the other hand,
setting $\delta_E\gg1$ leads C-GLISp to explore the entire domain $\domain$ regardless of
the decision-maker's preferences and feasibility/satisfaction assessments.

Functions $\surr{G}$ and $\surr{S}$ in~\eqref{eq:acquisition_modified} aim at discouraging  exploration in regions  where the experiment is  predicted to be infeasible (\emph{i.e.}, $\parvec \notin \domainHARD$) and/or unsatisfactory (\emph{i.e.}, $\parvec \notin \domainSOFT$). Therefore, a poor selection of the hyperparameters $\delta_G$ and $\delta_S$ and/or a poor predictive capability of $\surr{G}$ and $\surr{S}$ (\emph{e.g.}, due to a limited number of samples) can prevent finding
new vectors that are actually feasible and/or satisfactory. To alleviate this issue, we suggest to adaptively tune $\delta_G$ and $\delta_S$ based on the sampled standard deviation  obtained from leave-one-out cross-validation~\cite{stone1974cross} of $\surr{G}$ and $\surr{S}$, respectively. More specifically, each available sample in the set $\{\parvec_1,\ldots,\parvec_N\}$ is used once as a testing point and the remaining ones are used to train  $\surr{G}$ and $\surr{S}$. The prediction $\surr{G}(x_i)$ and $\surr{S}(x_i)$ on the test sample $x_i$ is compared with the corresponding  labels $G(x_i)$ and $S(x_i)$ assigned by the decision-maker to compute the following sampled standard deviations of the error:
\begin{align}
	\hat \sigma_G =  \min\left\{1,\ \sqrt{ \frac{ \sum_{i=1}^N (\surr{G}(x_i) - G(x_i))^2}{N-1} }\ \right\}, \nonumber \\
	\hat \sigma_S =  \min\left\{1,\ \sqrt{ \frac{ \sum_{i=1}^N (\surr{S}(x_i) - S(x_i))^2}{N-1} }\ \right\}.	
\label{eq:delta_tune-1}
\end{align}

The sampled standard deviations are then used to update, after each iteration,  the weights $\delta_G$ and $\delta_S$  as follows:
\begin{subequations} \label{eq:delta_tune-2}
\begin{align}
	\delta_{G} = & \ (1-\hat \sigma_G)\delta_{G,\text{default}},  \\
	\delta_{S} = & \ (1-\hat \sigma_S)\delta_{S,\text{default}},
\end{align}
\end{subequations}
where $\delta_{G,\text{default}}$ and $\delta_{S,\text{default}}$ are default values set by the user. Clearly, one should select $\delta_{G,\text{default}} > \delta_{S,\text{default}}$, so that infeasibility is penalized more than unsatisfactory behavior. The updated values of $\delta_G$ and $\delta_S$ are then used to construct the acquisition function $a(\parvec)$ in~\eqref{eq:acquisition_modified}.

The next sample $\parvec_{N+1}$ to test is obtained  by minimizing $a(\parvec)$, \emph{i.e.}, 
\begin{equation}
\parvec_{N+1}=\arg\min_{\parvec\in \domain%g(x)\leq 0
} a(\parvec).
\label{eq:xNp1}
\end{equation}
Different optimization methods can be used to solve problem~\eqref{eq:xNp1} efficiently either via derivative-free~\cite{RS13}, or derivative based algorithms.

C-GLISp updates the surrogates $\hat f$, $\surr{G}$, and $\surr{S}$, and the exploration function $z_N(\parvec)$, by iteratively suggesting a new point $\parvec_{N+1}$ to test, and by receiving feedback from the decision-maker in terms of feasibility,  overall satisfaction,  and preferences between pairs of experiments. % and    algorithm can approach optimal solution as number of experiments increases. 
Algorithm~\ref{algo:idwgopt-pref} summarizes the proposed method.

\begin{algorithm}[ht!]
	\caption{C-GLISp: Preference learning algorithm with unknown constraint handling }
	\label{algo:idwgopt-pref}
	~~\textbf{Input}: Lower and upper bounds $(\ell,u)$, known constraint set if available;
	number $N_{\rm init}\geq 2$ of initial samples, number $N_{\rm max}\geq N_{\rm init}$ of maximum function 
	evaluations; $\delta_E\geq 0$, $\delta_{G,\text{default}}\geq 0$ and $\delta_{S,\text{default}}\geq 0$; $\sigma>0$; and $\epsilon > 0$; 
	self-calibration index set  $\mathcal{I}_{\rm{sc}} \subseteq \{1,\ldots,N_{\rm max}-1\}$.
	\vspace*{.1cm}\hrule\vspace*{.1cm}
	\begin{enumerate}%[label*=\arabic*., ref=\theenumi{}]
		
		\item \label{algo:latin} Generate $N_{\rm init}$ random samples $X=\{x_1,\ldots,x_{N_{\rm init}}\}$ using Latin hypercube sampling method~\cite{MBC79};
				
		\item  $N\leftarrow 1$,  $i^\star\leftarrow 1$;
		\item \textbf{While} $N< N_{\rm max}$ \textbf{do}
		\begin{enumerate}%[label=\theenumi{}.\arabic*., ref=\theenumi{}.\arabic*]
		\item \textbf{if} $N =1$ \textbf{then}
			\begin{enumerate}%[label=\theenumi{}.\arabic*., ref=\theenumi{}.\arabic*]
				\item Observe feasibility $G_N$ and satisfaction $S_N$;
			\end{enumerate}
			\item \textbf{if} $N \geq N_{\rm init}$ \textbf{then}
					\begin{enumerate}%[label=\theenumii{}.\arabic*., ref=\theenumii{}.\arabic*]
					\item \label{algo:AL1} \textbf{if} $N \in \mathcal{I}_{\rm sc}$ \textbf{then} recalibrate $\epsilon$ through $K$-fold cross-validation;	
					\item  Solve~\eqref{eq:QP-pref} to obtain $\beta$ to define the surrogate function $\hat f$~\eqref{eq:rbf};
					\item Update $\delta_G$ and $\delta_S$ as in~\eqref{eq:delta_tune-2};
					\item \label{algo:acquisition} Define  acquisition function $a$ as in~\eqref{eq:acquisition_modified};
					\item \label{algo:globopt} Solve  optimization problem~\eqref{eq:xNp1} and get $x_{N+1}$;	
					\end{enumerate}
			\item $i(N)\leftarrow i^\star$, $ j(N) \leftarrow  N+1$;
			\item Observe feasibility $G_{j(N)}$. satisfaction $S_{j(N)}$ and preference $b_N = \pref(x_{i(N)},x_{j(N)})$ ;
			%\item Observe satisfaction $S_N$;
			%\item Observe preference $b_N = \pref(x_{i(N)},x_{j(N)})$; 
			\item \textbf{if} $b_N=1$ \textbf{then set} $i^\star\leftarrow  j(N)$; 
			\item $N\leftarrow N+1$;
		\end{enumerate}
		\item \textbf{End}.
	\end{enumerate}
	\vspace*{.1cm}\hrule\vspace*{.1cm}
	~~\textbf{Output}: Computed best input   $x^\star=x_{i^\star}$.
\end{algorithm}

%%%%%%%%%%%%%%%%%%%%%%%%%%%%%%%%%%%%%%%%%%%%%%%%%%%%%%%%%%%

\section{Optimization Benchmarks}\label{sec:benchmarks}
\begin{table*}[!bt]
	\caption{Numerical benchmarks  - Problem Specification}\label{tab:Benchmarks_spec}
	\centering
	\begin{adjustbox}{width=1\textwidth}
%	\vspace{-0.2cm}
	\begin{tabular}{lllll}
		\hline
		Test function & Objective function  & Unknown constraints & Search domain $\domain$ \\
		\hline
		Mishra's Bird function-     & $f(x,y) = \sin(y)e^{(1-\cos(x))^2} + \cos(x)e^{(1-\sin(y))^2}$         & Feasibility constraints:  & $[-10.0,-2]; $  \\ 
		    constrained (modified) \cite{mishra2006some,Bird_phenix} &  $\ \ \ \ \ \ \ \ \ \ \ \ \ \ \ \ \ \ \ \ \ \ \ \ \ \ \ \ \ \ \ \ \ \ \ \ \ \ +(x-y)^2$ & $(x+9)^2 + (y+3)^2<9$ & $[-6.5,0.0]$ \\
		   (MBC) & & & \\
		\hline
		camelsixhumps-   &  $f(x,y)=(4-2.1x^2+x^{4/3})x^2$    & Feasibility constraints:  $g_1  \cap   g_2$ & $[-2,2]; $      \\
		hard constrained \cite{Jamil2013-tr,Bem20}  & $\ \ \ \ \ \ \ \ \ \ \ \ \ \ \ \ \ \ \ \ \ \ \ \ \ \ \ \ \ \ \ \ \ \ \ \ \ \ +xy+(4y^2-4)y^2$ &  & $[-1,1] $  \\ 
		
 (CHC) & &$g_1:\left[\begin{smallmatrix}1.6295 & 1\\ -1 & 4.4553\\ -4.3023 & -1\\ -5.6905 & -12.1374\\ 17.6198 & 1 \end{smallmatrix}\right]\left[\begin{smallmatrix} x\\ y \end{smallmatrix}\right] <\left[\begin{smallmatrix}3.0786 \\ 2.7417 \\ -1.4909 \\ 1 \\ 32.5198 \end{smallmatrix}\right] $ &  \\ 
  & & $g_2: x^2 + (y+0.1)^2 < 0.5$ &  \\
		\hline
		camelsixhumps-        &$f(x,y)=(4-2.1x^2+x^{4/3})x^2$        & Feasibility constraints: $g_2$& $[-2,2];$  \\
		hard and soft constrained \cite{Jamil2013-tr,Bem20}& $\ \ \ \ \ \ \ \ \ \ \ \ \ \ \ \ \ \ \ \ \ \ \ \ \ \ \ \ \ \ \ \ \ \ \ \ \ \ +xy+(4y^2-4)y^2$  & Satisfaction constraints: $g_1$ & $[-1,1] $  \\
		 
(CHSC) & &$g_1: \left[\begin{smallmatrix}1.6295 & 1\\ 0.5 & 3.875\\ -4.3023 & -4\\ -2 & 1\\ 0.5 & -1 \end{smallmatrix}\right] \left[\begin{smallmatrix} x\\ y \end{smallmatrix}\right] <\left[\begin{smallmatrix}3.0786 \\ 3.324 \\ -1.4909 \\ 0.5 \\ 0.5 \end{smallmatrix}\right] $ &  \\ 
  & &$g_2:  x^2 + (y+0.04)^2 < 0.8$ &  \\

		\hline
	\end{tabular}
	\end{adjustbox}
%	\vspace{-0.4cm}
\end{table*}

\begin{table*}[!bt]
	\caption{Numerical benchmarks - Solver Specification}\label{tab:Benchmarks_para}
	\centering
	\begin{adjustbox}{width=1\textwidth}
%	\vspace{-0.2cm}
	\begin{tabular}{lcccccccccccc}
		\hline
		Test function & Max number of & Number of initial & \multicolumn{3}{c}{Hyper-parameter values}  & & \multicolumn{3}{c}{RBF specifications~\eqref{eq:rbf}}& Tolerance  & Weights & Regularization\\
		\cline{4-6}
		\cline{8-10}
		& fun. eval. $N_{max}$ & sampling $N_{init}$ &  $\delta_E$& $\delta_{G,\text{default}}$  & $\delta_{S,\text{default}}$  && function & initial $\epsilon$  & recalibration steps & $\sigma$ in~\eqref{eq:QP-pref}& $c_h$ in~\eqref{eq:QP-pref} & $\lambda$ in~\eqref{eq:QP-pref}  \\
		\hline
		MBC     & 50 &13& 1.0 & 1.0& $-$& & Inverse quadratic& 1.0&$\{13,22,32,41\}$ &0.02 &1.0&1$e$-6 \\
		CHC & 100&25 & 2.0& 2.0& $-$ & & Inverse quadratic& 1.0& $\{25,44,63,81\}$ &0.01 &1.0&1$e$-6   \\
		CHSC  & 50 &13& 1.0& 1.0& 0.5& & Inverse quadratic& 1.0& $\{13,22,32,41\}$ &0.02&1.0&1$e$-6    \\

		\hline
	\multicolumn{13}{l}{Same parameters (if relevant) are used in C-GLISp, GLISp, and PBO.}\\
%	\multicolumn{12}{l}{b- Only $\delta$ is included in the formulation of GLISp algorithm. }\\
%	\multicolumn{10}{l}{c- Value in () indicates the number of tests the optimum satisfies the soft constraints}\\
	\end{tabular}
	\end{adjustbox}
%	\vspace{-0.4cm}
\end{table*}

\begin{table*}[!bt]
	\caption{Numerical benchmarks - Results}\label{tab:Benchmarks_result}
	\centering
	\begin{adjustbox}{width=1\textwidth}
%	\vspace{-0.2cm}
	\begin{tabular}{p{0.12\textwidth}*{4}{>{\centering}p{0.1\textwidth}}p{0.005\textwidth}*{4}{>{\centering}p{0.1\textwidth}}}
		\hline
		Test function & \multicolumn{4}{c}{Constrained optimum$^{a}$} & &\multicolumn{3}{c}{Feasibility$^{b}$}
		\tabularnewline
		\cline{2-5}
		\cline{7-9}
		& Optimum&  C-GLISp  & GLISp & PBO & &  C-GLISp & GLISp & PBO 
		\tabularnewline
		\hline
		MBC    & -48.4& -47.95  & -48.33  & -40.24  & & 100 & 100 & 91  
		\tabularnewline
		CHC & -0.5844   &-0.3582  & -0.5224 & 0.2571 & & 96 & 66 & 33 
		\tabularnewline
		CHSC  & -0.9050&  -0.8526   & -0.8861 & -0.6315  & & 96 (95) & 82 (84) & 74 (72)
		\tabularnewline

		\hline
	\multicolumn{9}{l}{a- The median of computed constrained optima that are feasible out of $100$ runs (the distribution over 100 runs is reported in  Table~\ref{tab:100_tests_result_dist}).}\\
	\multicolumn{9}{l}{b- Number of runs whose computed optimizers are feasible out of 100 runs. Values in parentheses indicate the number of runs the optimizer is satisfactory.} \\
	\end{tabular}
	\end{adjustbox}
\end{table*}

We test C-GLISp on three constrained global optimization benchmarks to illustrate its effectiveness in solving optimization problems with unknown constraints. Computations are performed on an Intel i7-8550U 1.8-GHz CPU laptop with 8GB of RAM. The  Latin hypercube sampling method~\cite{MBC79} (\emph{lhsdesign} function of  the  \emph{Statistics  and  Machine  Learning Toolbox} of MATLAB~\cite{matlabStaML}) is used in the initial sampling phase of C-GLISp. Particle Swarm Optimization (PSO)~\cite{VV09} is used to minimize the acquisition function 
as in~\eqref{eq:xNp1}.

C-GLISp is compared to the original GLISp and to PBO (with \emph{expected improvement} as acquisition function)~\cite[Section 2.3]{brochu2007active}. For numerical benchmarks, C-GLISp, GLISp, and PBO assign the preferences on pairwise comparisons based on the combined assessments of the objective function value, feasibility, and performance satisfaction. For each test function, depending on the problem formulation, a maximum of three types of queries is obtained when using C-GLISp, which are the preference relation ($B$), the feasibility label ($G_F$), and the satisfaction label ($S_F$). In contrast, GLISp and PBO only rely on the preference relation $B$. The goal of the comparison between C-GLISp and GLISp is to check if accounting the feasibility and/or satisfactory information explicitly in the acquisition function can encourage safe exploration from the comparisons. It is worth noting that the exact evaluation of the objective function and the constraints (feasibility and satisfactory outcomes) for these numerical benchmarks are unknown to the algorithms and are only used to construct a synthetic decision-maker. 

Table~\ref{tab:Benchmarks_spec} lists the specifications of the benchmarks. The original feasibility set  of the test function \emph{Mishra's Bird function-constrained} (MBC)~\cite{mishra2006some,Bird_phenix} is modified so that the unconstrained global optimum in the search domain is no longer in the feasible area. The \emph{camelsixhumps-hard constrained} (CHC) benchmark~\cite{Jamil2013-tr,Bem20} considers two feasibility constraints, and the unconstrained global optimum also differs from the constrained one. Lastly, the benchmark function \emph{camelsixhumps-hard and soft constrained} (CHSC)~\cite{Jamil2013-tr,Bem20} has both feasibility and satisfaction constraints. The two unconstrained optima for this test function are both feasible but not satisfactory.

The values of the hyper-parameters in C-GLISp, GLISp, and PBO are provided in Table~\ref{tab:Benchmarks_para}. The number of initial samples ($N_{init}$) is selected as one fourth of the maximum number of function evaluations ($N_{max}/4$) rounded to the nearest integer. Three-fold cross-validation is used to update the hyper-parameter $\epsilon$~\eqref{eq:rbf} at iterations $N_{init}$, $N_{init}+(N_{max}-N_{init})/4$, $N_{init}+(N_{max}-N_{init})/2$, and $N_{init}+3(N_{max}-N_{init})/4$, rounded to the closest integers, which
define the self-calibration index set  $\mathcal{I}_{\rm{sc}}$. The tolerance $\sigma$ in~\eqref{eq:QP-pref} is set to $1/N_{max}$. The default value $\delta_{G,\text{default}}$ in~\eqref{eq:delta_tune-2} is the same as $\delta_E$, so that the feasibility term in~\eqref{eq:acquisition_modified} is comparable to the pure exploration term, while the default value  $\delta_{S,\text{default}}$ is selected as $\delta_{G,\text{default}}/2$ to reduce its effects with respect  to hard feasibility constraints. The parameters $\delta_G$ and $\delta_S$ are kept at their default values during the first $N_{init}$ experiments, then updated each time a new point is added using equation~\eqref{eq:delta_tune-2}. The remaining parameters of the solvers are set according to the defaults used or suggested in~\cite{BemPig20}. 

\begin{figure}
  \centering
  \includegraphics[width=\linewidth]{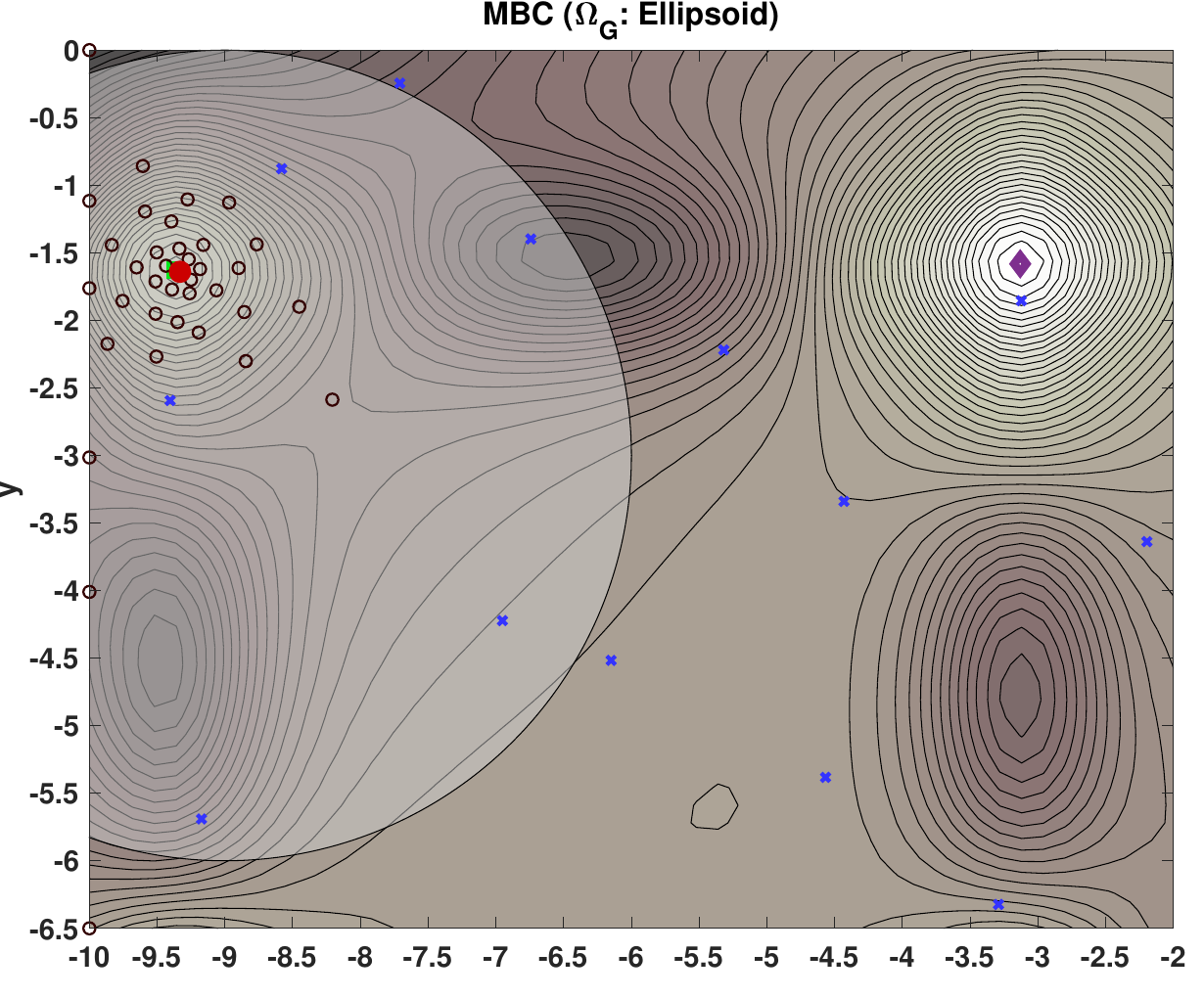}
  \includegraphics[width=\linewidth]{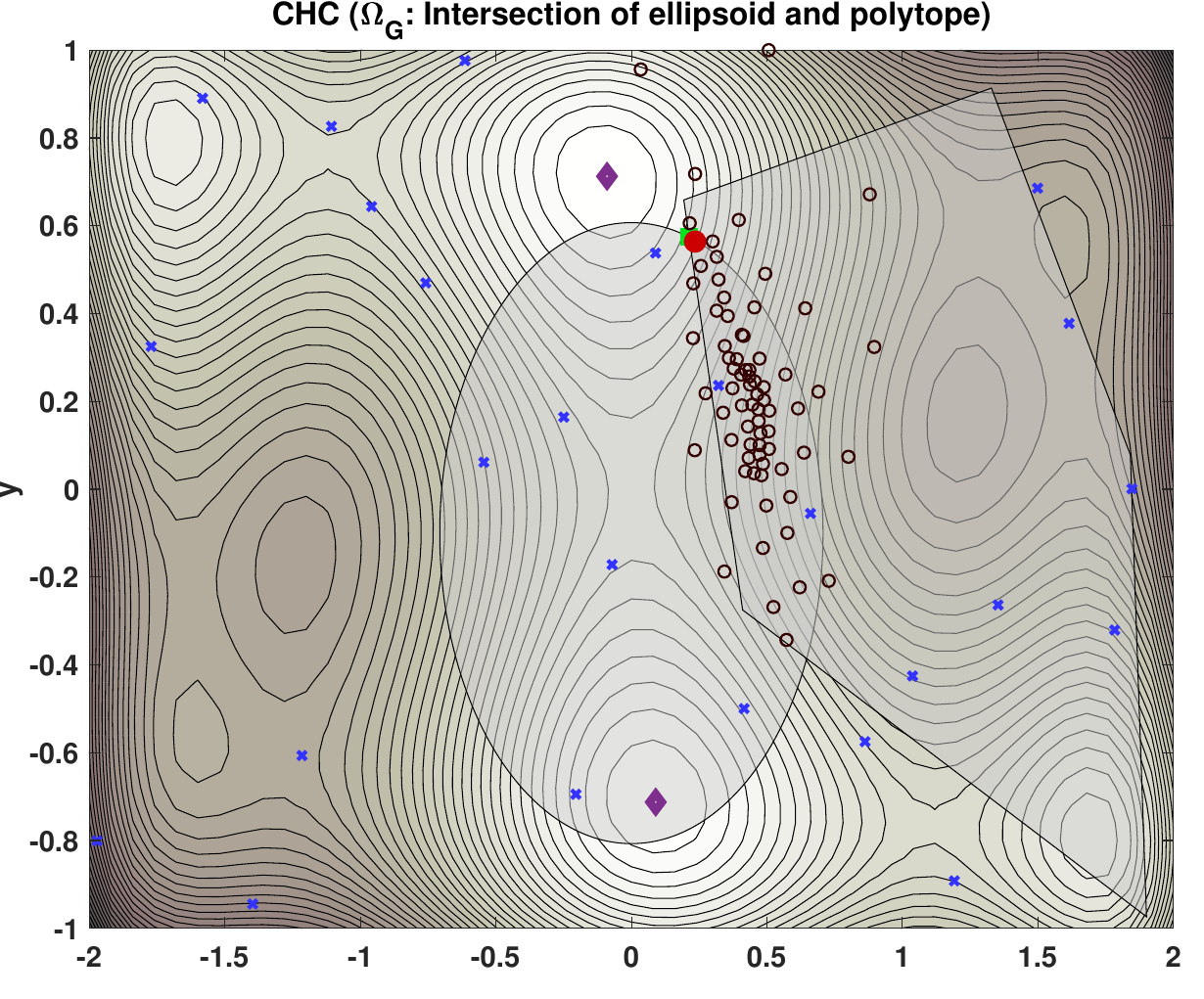}  
  \includegraphics[width=\linewidth]{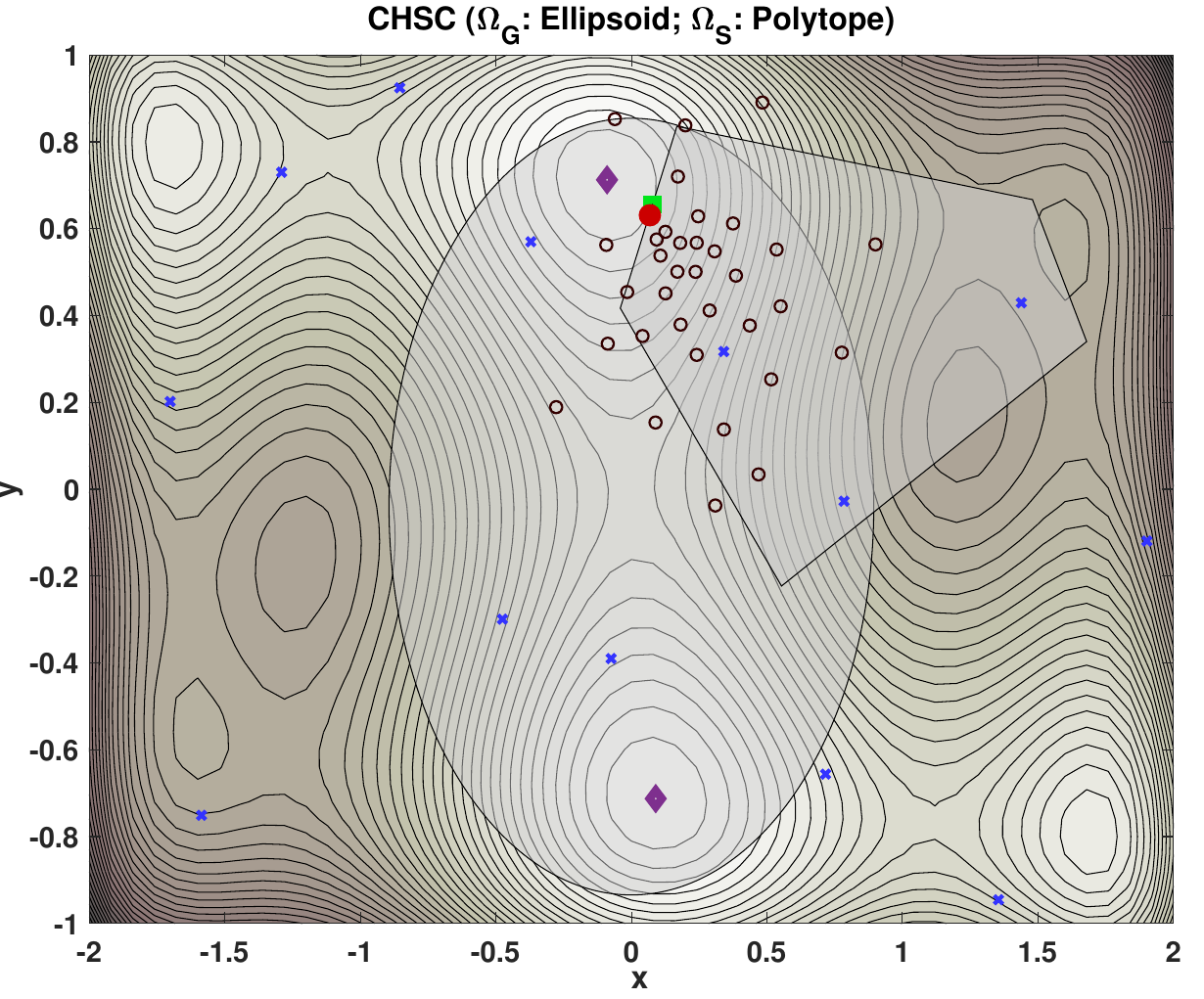}
  
%  \vspace{-0.7cm}
  \caption{Algorithm C-GLISp. Level sets of the functions used in the three benchmarks, along with feasibility and satisfaction sets.  Blue $ \times$: points generated from initial sampling phase; black $\circ$: points generated from active learning phase;  purple $\Diamond$: global unconstrained optimizer; red $\bullet$: constrained optimizer found after $N_{max}$ iterations; green $\square$: global constrained optimizer. As $N$ increases, the points generated by C-GLISp approach the constrained optimizer, and most of the points generated during the active learning phase lay in the feasibility and satisfaction regions.} 
  \label{fig:benchmark_one_test}
  %\vspace{-0.4cm}
\end{figure}

Table~\ref{tab:Benchmarks_result} reports the results obtained by running a Monte-Carlo simulation with 100 runs of C-GLISp, GLISp, and PBO to obtain statistically significant results. One of such runs of C-GLISp on all three numerical benchmarks is depicted in Fig.~\ref{fig:benchmark_one_test}.  Table~\ref{tab:100_tests_result_dist} displays the distribution over 100 runs of the percentage difference between the achieved feasible solutions and the true constrained optimum.
Overall, the results from Table~\ref{tab:Benchmarks_result} and~\ref{tab:100_tests_result_dist} show that C-GLISp can find a feasible near-optimal solution more frequently than GLISp and PBO. 

\begin{table}[]
\caption{Distribution over 100 runs of the percentage difference between achieved and global optimum}\label{tab:100_tests_result_dist}
\centering
\begin{adjustbox}{width=0.47\textwidth}
\begin{tabular}{p{0.07\textwidth}p{0.07\textwidth}*{4}{p{0.06\textwidth}}}

\hline
\multirow{3}{*}{Benchmark} & \multirow{3}{*}{Algorithm} & \multicolumn{4}{c}{Number of runs within each interval}                  \\ \cline{3-6} 
                           &                            & \multicolumn{4}{c}{Intervals of \% Difference from Global Optimum} \\
                           &                            & (0,5{]}             & (5,10{]}             & (10,15{]}            & (15,100{]}            \\ \hline
\multirow{3}{*}{MBC}       & PBO                        & 39                  & 4                    & 2                    & 18                    \\
                           & GLISp                      & 67                  & 1                    & 2                    & 3                     \\
                           & C-GLISp                    & 69                  & 6                    & 1                    & 5                     \\ \hline
                           &                            & (0,5{]}             & (5,20{]}             & (20,50{]}            & (50,100{]}            \\ \hline
\multirow{3}{*}{CHC}       & PBO                        & 0                   & 0                    & 4                    & 7                     \\
                           & GLISp                      & 28                  & 7                    & 5                    & 1                     \\
                           & C-GLISp                    & 0                   & 20                   & 40                   & 22                    \\ \hline
                           &                            & (0,5{]}             & (5,10{]}             & (10,15{]}            & (15,100{]}            \\ \hline
\multirow{3}{*}{CHSC}      & PBO                        & 13                  & 10                   & 4                    & 27                    \\
                           & GLISp                      & 56                  & 8                    & 5                    & 9                     \\
                           & C-GLISp                    & 43                  & 22                   & 13                   & 16                    \\ \hline
\multicolumn{6}{l}{Note: only the runs with feasible solutions within 100\% difference from the   }                                                                                                                           \\
\multicolumn{6}{l}{global optimum are counted.}  
\end{tabular}
\end{adjustbox}
\end{table}

\begin{figure}
  \centering
  \includegraphics[width=\linewidth]{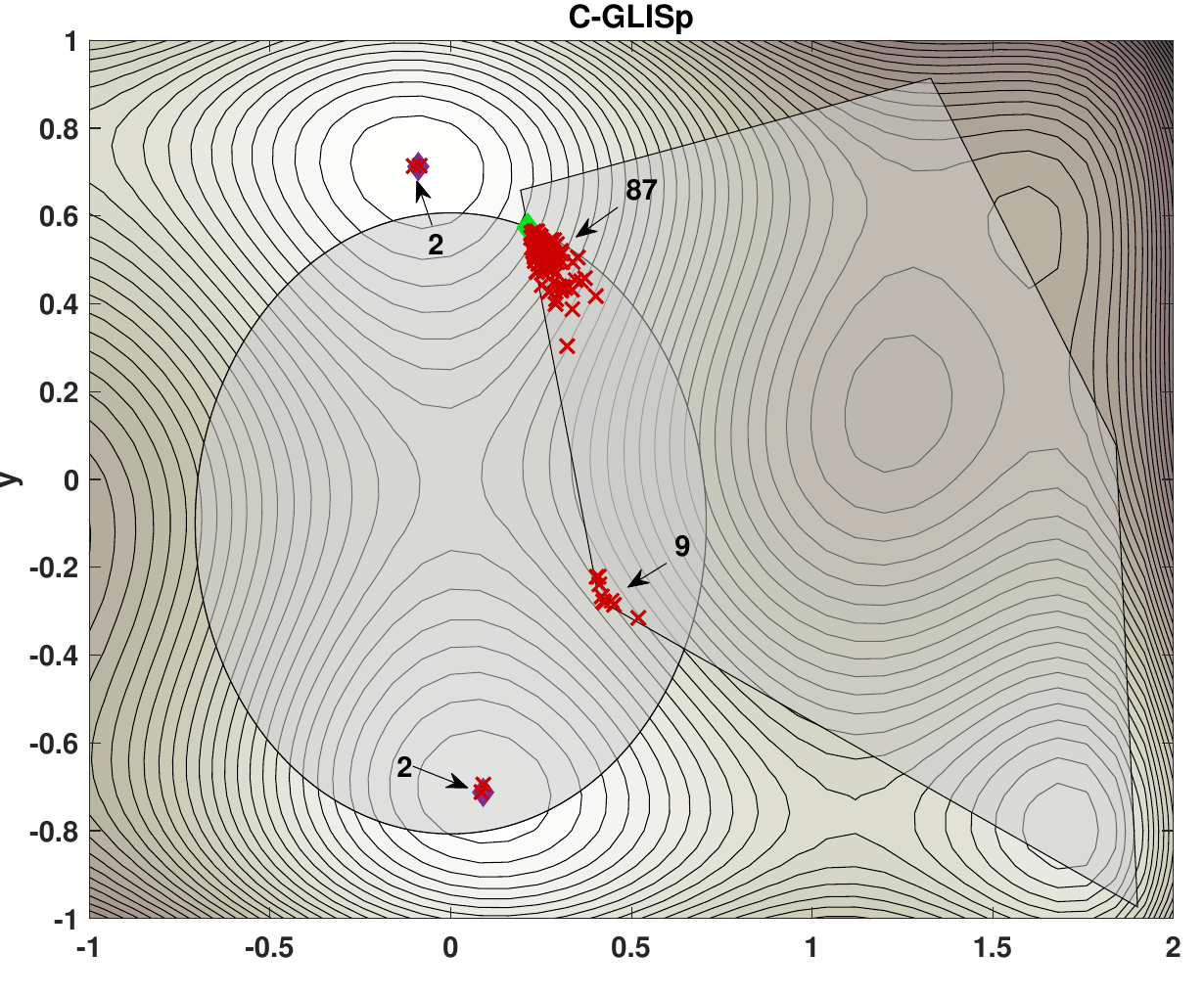}
  \includegraphics[width=\linewidth]{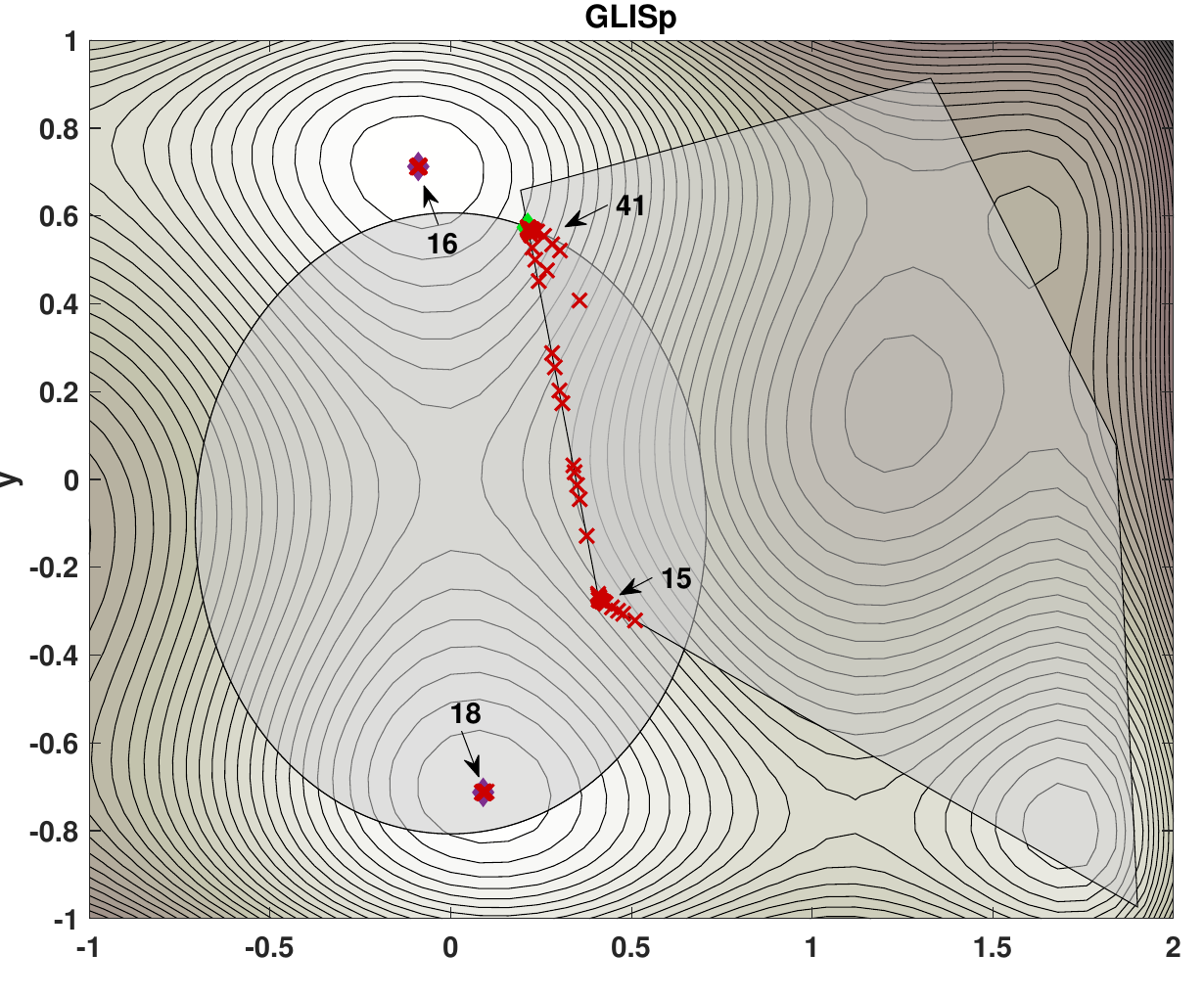}
    \includegraphics[width=\linewidth]{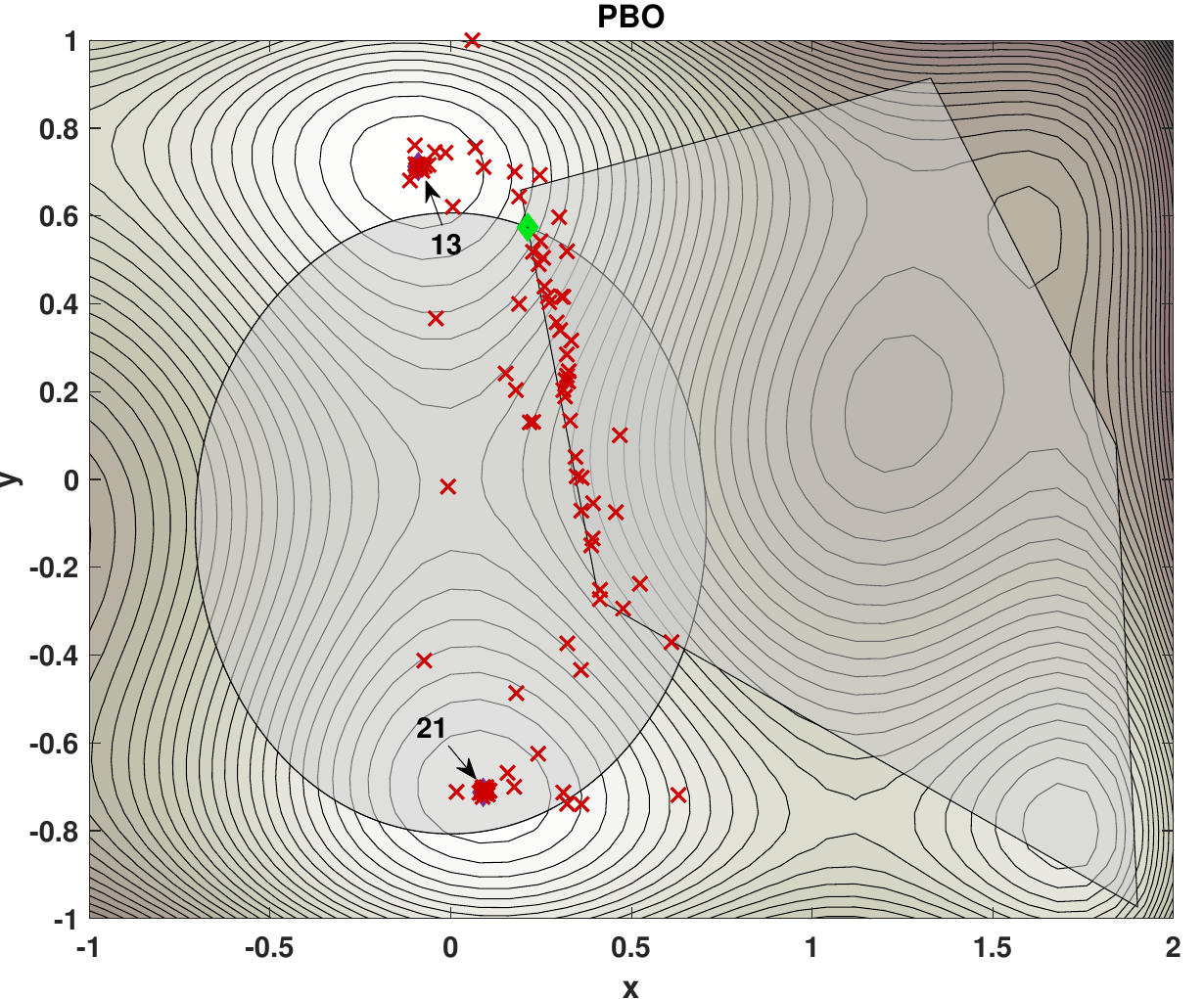}
%  \vspace{-0.7cm}
  \caption{Benchmark CHC. Optimizers computed  by C-GLISp, GLISp and PBO in 100 runs. Red $\times$: optimizer computed at the end of each run; purple $\Diamond$:  unconstrained optimizer; green $\Diamond$: global constrained optimizer. Numbers in black with arrows indicate the number of overlapping points. }
  \label{fig:CHC_100_tests}
\end{figure}
 
 From the results on  the benchmark function MBC, where the feasibility set $\domainHARD$  covers roughly one-third of the domain $\domain$ (cf. Fig.~\ref{fig:benchmark_one_test}), the performance of GLISp and C-GLISp are comparable. They always terminate the search with a feasible optimum (100 out of 100 runs), with 67$\%$ and 69$\%$ of them, respectively, located within 5$\%$ difference from the global solution (Table~\ref{tab:100_tests_result_dist}). 
On the other hand, PBO computes a feasible optimum  in 91 runs, but with only around 39$\%$ within 5$\%$ difference (Table~\ref{tab:100_tests_result_dist}). When the constraint is more complex such as the one in CHC, the majority of the optima computed by C-GLISp (96 out of 100 runs) are feasible. In comparison, only 66 and 33 runs by GLISp and PBO, respectively,  terminate with a feasible solution (Table~\ref{tab:Benchmarks_result} and Fig.~\ref{fig:CHC_100_tests}). From Fig.~\ref{fig:benchmark_one_test}, it is also observed that, for the test function CHC, after the initial sampling phase, most points generated in the active learning phase by C-GLISp are within the feasible region.

For the test function CHSC, C-GLISp often find a near-optimal solution that is both feasible and satisfactory. The performance of GLISp is slightly worse than C-GLISp in terms of the number of times a feasible and satisfactory solution is obtained. PBO can identify a feasible and satisfactory solution with a relatively high chance but still lower than both C-GLISp and GLISp. Also, its final outcome is worse, see Table~\ref{tab:Benchmarks_result}. 

Table~\ref{tab:Benchmarks_result} also shows that, within the same number of iterations, the median of the computed feasible constrained optima from GLISp is always closer to the global constrained optimum than the one computed from C-GLISp. This is because of the trade-off between trying to get a more accurate solution (which is often achieved by sampling multiple points close to the current best solution up to iteration $N$) and exploring a larger area to reduce uncertainty (in problems with unknown constraints, C-GLISp also tries to identify possible feasible regions). 

For our problem setting, we set a limit on the computational budget. Modification of the exploration term from~\eqref{eq:IDW-distance} to~\eqref{eq:z_comb} helps to better escape from local minima in the early iterations by encouraging the exploration of regions of $\domain$ further away from the current best solution. 
This modification is significant for problems with small feasible regions (relative to the search domain) and/or complex unknown constraints (\emph{e.g.}, numerical benchmark CHC). 
This is because the additional exploration introduced by the modification can help the solver identify the feasible region more quickly and start recommending feasible guesses faster, reducing the chance that the solver gets trapped into an infeasible local optimum.
 
From the number of feasible solutions computed shown in Table~\ref{tab:Benchmarks_result} and the computed optimizers displayed in Fig.~\ref{fig:CHC_100_tests}, we observe that the situation of trapping into an infeasible local optimum occurs to GLISp more frequently than to C-GLISp. However, GLISp  can achieve a solution closer to the constrained optimum than C-GLISp (Table~\ref{tab:100_tests_result_dist}) when GLISp successfully identifies the feasible region. This is because  more computational budget is then used to get closer to the optimum than to explore other potentially feasible regions as in C-GLISp.  For problems where testing an infeasible solution is expensive and/or dangerous, it is better to be conservative and have a less optimal but feasible solution. As a result, C-GLISp is preferred over GLISp under these problem settings.

Overall, the results on the numerical benchmarks show that both C-GLISp and GLISp can approach near-optimal solutions within a small number of function evaluations. In all of the three benchmarks, both  C-GLISp and GLISp outperform PBO (Table~\ref{tab:Benchmarks_result} and~\ref{tab:100_tests_result_dist}). An explanation for the superior performance of C-GLISp in identifying feasible/satisfactory solutions
is that it explicitly leverages feasibility/satisfaction information in the acquisition function,
while GLISp and PBO handle unknown constraints only through preference queries.

\section{MPC Calibration}\label{sec:case_study}
To illustrate the application of C-GLISp to controller calibration, we consider the design of an MPC
controller for lane-keeping (LK) and obstacle-avoidance (OA) in autonomous driving.  
MPC is employed to command vehicle velocity and steering angle to provide a smooth and safe drive.
The same problem was considered  in~\cite{zhu2020pref} and is extended in this paper to handle feasibility and satisfaction constraints. 
  
The design of model predictive controllers  requires tuning several knobs, such as the prediction and control horizons, the weight matrices in the cost function, numerical tolerances in the optimization solver, \emph{etc.}, under hard constraints such as finding solutions within the sample interval of a real-time implementation.  Thus, it is hard to well define in advance a single quantitative performance index that captures a usually multifaceted desired closed-loop behavior to tune the MPC parameters automatically. 
 
We use C-GLISp to implement an iterative semi-automated calibration procedure (automatic selection of the MPC parameters and human-based qualitative assessment of closed-loop performance), thus avoiding  the burden of defining a quantitative index to minimize. The calibrator only needs to express a \emph{preference} on pairwise comparisons to indicate which of the two MPC tunings  has better closed-loop performance. Different from the case study in~\cite{zhu2020pref}, where we account for the information of feasibility and satisfaction conditions implicitly in the preference query (GLISp) as if they were, implicitly, underlying penalty functions, C-GLISp  explicitly takes into account this information via direct queries as in~\eqref{eq:fes_fun} and~\eqref{eq:softcost_fun}.

\subsection{System Description} 
We consider a simplified two-degree-of-freedom bicycle model with the front wheel as the reference point to describe the vehicle dynamics and simulate the experiment. The state variables $s = [x_f \  w_f \  \theta]'$ in the model are the longitudinal $x_f$ and lateral  $w_f$ [m] positions of the front wheel, and the yaw angle $\theta$ [rad]. The manipulated variables $u = [v \ \psi]'$ are the commanded vehicle velocity $v$ [m/s] and steering angle $\psi$ [rad]. The standard kinematic equations

\begin{equation}\label{eq:nonlin_bicyclemodel}
\begin{split}
 \dot{x}_f = & v \cos(\theta+\psi)\\
 \dot{w}_f = & v \sin(\theta+\psi)\\
 \dot{\theta} = & \frac{v\sin(\psi)}{L}
\end{split}
\end{equation}
are used to model the evolution of the vehicle,
where $L$ [m] is the vehicle length. Here, full state observation is assumed, \emph{i.e.}, the control output $y = s $.

\subsection{MPC Formulation}\label{sec:mpc_desp}
The formulation of the semi-automatic MPC calibration process is detailed in~\cite{zhu2020pref}. The dynamical  model~\eqref{eq:nonlin_bicyclemodel} is linearized around its nominal point $\bar{s}_k = [\bar{x}_{f_k} \  \bar{w}_{f_k} \  \bar{\theta}_k]' $, $\bar{u}_k = [\bar{v}_k \ \bar{\psi}_k]'$, and $\bar{y}_k = \bar{s}_k$  at each time step and discretized with sampling time $T_s$, resulting the following discrete-time state-space model:
\begin{equation}  
	%\label{eq:SITO}
	\begin{split}
	 \tilde{s}_{k+1}&\! = \! \left[\begin{smallmatrix} 1 & 0 & - \bar{v}_k \sin(\bar{\theta}_k+\bar{\psi}_k) T_s \\0 & 1 & \bar{v}_k \cos(\bar{\theta}_k+\bar{\psi}_k) T_s\\0 & 0 & 1 \end{smallmatrix}\right]\tilde{s}_k\!\\
	 &\quad +\!  
	 \left[\begin{smallmatrix} \cos(\bar{\theta}_k+\bar{\psi}_k) T_s & - \bar{v}_k \sin(\bar{\theta}_k+\bar{\psi}_k) T_s \\ \sin(\bar{\theta}_k+\bar{\psi}_k) T_s & \bar{v}_k \cos(\bar{\theta}_k+\bar{\psi}_k) T_s \\ \frac{\sin(\bar{\psi}_k)}{L}T_s & \frac{\bar{v}_k\cos(\bar{\psi}_k)}{L}T_s  \end{smallmatrix}\right] \tilde{u}_k  \\ %\label{eq:MPCsysstate}
	 \tilde{y}_k& \! =\tilde{s}_k,   
	\label{eq:MPCsys}
	\end{split}
\end{equation}
%\end{minipage}}
where subscript $k$ denotes the value at time step $k$ and $\widetilde{\text{Var}} = \text{Var} - \overline{\text{Var}}$. This prediction model is then used to design a linear MPC via a real-time iteration scheme~\cite{GZQBD16,diehl2005realtimeIter}. At each sampling time $t$, the following QP problem is solved to compute the MPC action to be applied:
\begin{equation} \label{eq:MPC}
	\begin{split}
	  \min_{\left\{u_{t+k|t}\right\}_{k=0}^{\Nu-1}}
 %, \varepsilon}  
	\sum_{k=0}^{\Np-1} (\left\|y_{t+k|t}-y^{\rm{ref}}_{t+k}\right\|_ {Q_y}^2 + \left\| \Delta u_{t+k|t} \right\|_{Q_{\Delta u}}^2) \\ 
	%\label{eq:MPCcost}
	%& \qquad   \qquad+\! \! \sum_{k=0}^{\Np-1} \left\| \Delta u_{t+k|t} \right\|_{Q_{\Delta u}}^2 \\ 
	%& \qquad   \qquad+\! \! \sum_{k=0}^{\Np-1} \left\|u_{t+k|t}-u^{\rm ref}_{t+k}\right\|_{Q_u}^2  \! + \! \sum_{k=0}^{\Np-1} \left\| \Delta u_{t+k|t} \right\|_{Q_{\Delta u}}^2 \\ 
% & \qquad  \qquad	\mathrm {s.t. } \qquad  \mathrm{input\ and\ output\ constraints} \\
%	 & \qquad  \qquad \qquad  \ \ \ \  \mathrm{model\ equation\   \eqref{eq:MPCsys}} 
%	%+ Q_\varepsilon \left\| \varepsilon\right\|^2 \\%+ Q_{\rK}
	\end{split}
\end{equation}
subject to model equation~\eqref{eq:MPCsys} and the following input and output constraints
\begin{equation}  \label{eq:MPC_constraints}
	\begin{split}
	& \quad y_{\mathrm{min}} \leq y_{t+k|t} \leq  y_{\mathrm{max}}, \  k=1,\ldots,\Np  \\
	& \quad   u_{\mathrm{min}}  \leq u_{t+k|t} \leq  u_{\mathrm{max}}, \  k=1,\ldots,\Np \\
	& \quad  \Delta u_{\mathrm{min}}   \leq \Delta u_{t+k|t}\leq  \Delta u_{\mathrm{max}}, \ \    k=1,\ldots,\Np   \\
	%\label{eq:MPC:Ducons1}
%	& \quad    \Delta u_{t+k|t} \leq  \Delta u_{\mathrm{max}}, \ \  k=1,\ldots,\Np \\ 
	%\label{eq:MPC:Ducons2} %\label{eq:MPC:ucons} \\
	& \quad  u_{t+\Nu+j|t}= u_{t+\Nu|t},  \ \  j=1,\ldots, \Np-\Nu, 
	%\label{eq:MPCconst}
	\end{split}
\end{equation}
where $\left\|\cdot\right\|^2_M$ is the squared norm weighted by a matrix $M$; $\Delta u_{t+k|t} =  u_{t+k|t}- u_{t+k-1|t}$; $y_{\rm ref}$ and $u_{\rm ref}$ are the reference values of control outputs and inputs during the experiment (which are unique to the LK and OA phases of the experiment); and $N_p$ and $N_u$ are the prediction and control horizons.

\subsection{Control Objectives}
The two main objectives involved in this control task are: (1) maintain the vehicle in the same lane with constant speed if no obstacles (other vehicles) are present; and (2) pass other moving vehicles if they are within a safety distance. The qualitative descriptions of these  objectives, such as the ambiguity of transferring ``optimal obstacle avoidance'' into a mathematical formula, make it challenging to define a proper quantitative performance index for closed-loop performance. On the other hand, it is easier for a calibrator to compare the outcome of two different driving tests and then express a preference.

We describe the test scenario as follows. Note, for ease of assessment, the unit of $v$ and $\psi$ described in the following text as well as in the figures are represented in km/hr and degree ($^\circ$), respectively. The controlled vehicle is initially at position ($x_f$, $w_f$) = (0, 0) m with $\theta$ = 0$^\circ$. Another vehicle (obstacle) is at position (30, 0) m and moving horizontally at a constant speed of 40 km/hr. The shape of both vehicles is assumed to be rectangular, with a length of 4.5 m and a width of 1.8 m. During nominal LK conditions, the vehicle being controlled moves horizontally at 50 km/hr, with $w_f$ = 0 m and $\dot{w}_f$ = 0 m/s . Once the obstacle is within a  safety distance, the vehicle being controlled should pass it while keeping a safe lateral distance between them. In this case, the horizontal and lateral safety distances are 10 m and 3 m, respectively. The vehicle controlled by MPC can vary its velocity in the range of [40, 70] km/hr during the LK period and [50, 70] km/hr during the OA period, and its reference velocities are set to $50$ and $60$ km/hr, respectively. For both LK and OA periods, $\theta$ can take  values in the range of [-45, 45]$^\circ$, with its rate of change between each time step limited to [-5, 5]$^\circ$/s. 

The following MPC design parameters are tuned:  sampling time ($T_s$);  prediction and control horizons ($N_p, N_u$);   weight matrix $Q_{\Delta u}$ in~\eqref{eq:MPC}. The sampling time is restricted in the interval [0.085, 0.5]~s. The prediction horizon is an integer allowed to vary in the range [10,30] and the control horizon is taken as a fraction $\epsilon_c$ of $N_p$ rounded up to the closest integer, with $\epsilon_c \in$ [0.1, 1]. The weight matrix is set to be diagonal $Q_{\Delta u} = \left[\begin{smallmatrix}q_{u11} & 0 \\0 & q_{u22} \end{smallmatrix}\right] $ and the values of  $\log(q_{u11})$ and $\log(q_{u22})$ are restricted in the interval [-5, 3]. The rest of the MPC design parameter is fixed, with $Q_y=\left[\begin{smallmatrix}1 & 0 & 0 \\0 & 1 & 0\\0 & 0 & 0 \end{smallmatrix}\right]$.

\subsection{Calibration Process}
The first author of the paper plays the role of the calibrator for this case study. The maximum number of function evaluations $N_{max}$ is set to 50,  with $N_{init}=$ 10.  The default hyperparameters $\delta_E$, $\delta_G$ and $\delta_S$ in~\eqref{eq:acquisition_modified} are set to 1, 1, and 0.5, respectively. The parameters $\sigma$, $c_h$ and $\lambda$ in~\eqref{eq:QP-pref} are set to 0.02, 1, and 1e-6, respectively. The hyperparameter $\epsilon$ characterizing  the RBF function~\eqref{eq:rbf} is initialized to 1.0, and  recalibrated at iterations 10, 20, 30, and 40 via 3-fold cross-validation. 

\begin{figure}[t!]
  \centering
  \includegraphics[width=\linewidth]{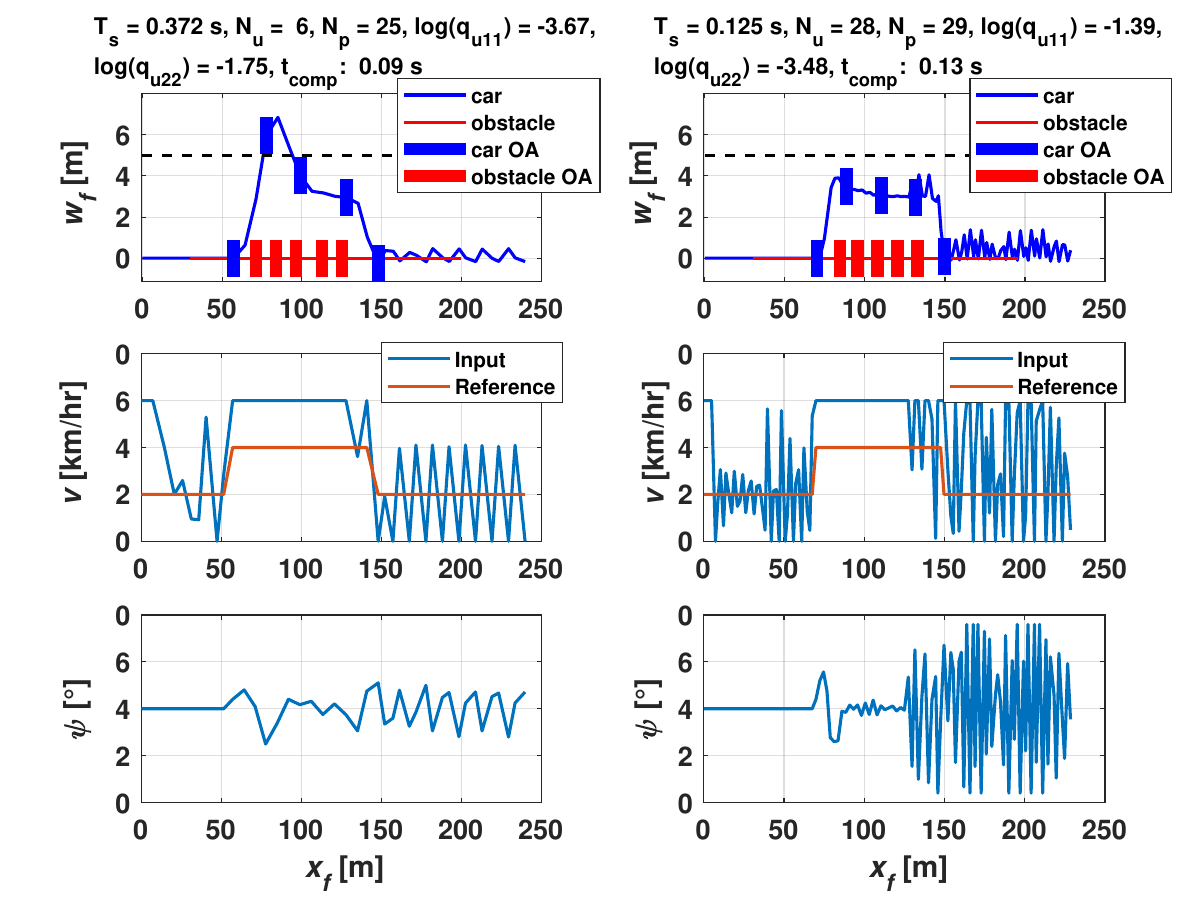}
  \vspace{-0.7cm}
  \caption{Vehicle control query window. The top subplots show the vehicle and obstacle positions. The ``vehicle OA'' and ``obstacle OA'' bars show five relative positions of the vehicle and obstacle during the OA period. The dashed lines indicate the lateral distance that the car should avoid exceeding (5 m in this case). The middle subplots show the actual and reference velocity $v$ at different longitudinal positions. The steering angle $\psi$ over the longitudinal position is depicted in the bottom subplots.
  	The results on the left panels are preferred and feasible, while the results on the right panels are infeasible. The results on both sets of panels fail to satisfy the satisfaction conditions. }
  \label{fig:car_comparison}
  %\vspace{-0.4cm}
\end{figure}

The closed-loop experiment of the vehicle control is simulated for 15 seconds. Fig. \ref{fig:car_comparison} shows the query window for one iteration of the calibration process. The MPC design parameters and the worst-case computational time ($t_{comp}$) required for solving the QP problem of MPC~\eqref{eq:MPC} are displayed at the top of the figure. At each iteration, the calibrator is asked to decide whether the newly proposed experiment is feasible and/or satisfactory, and to express a preference between the new experiment and the current best one. % achieved  select the preferred controller through pairwise comparisons 
More specifically, the calibrator labels the control policy that leads to ``unstable/unsafe'' and ``unimplementable'' behavior as infeasible (\emph{i.e.}, $G(x)=0$). Examples of ``unstable/unsafe'' behaviors include but are not limited to: vehicle hits the obstacle; vehicle oscillates on the road, \emph{etc.} ``Unimplementable'' cases are the ones whose computational time ($t_{comp})$ required for solving the QP problem of MPC~\eqref{eq:MPC} exceeds the sampling time $T_s$. As for being labeled as satisfactory, the following two criteria are used: ($i$) the lateral position of the vehicle does not exceed 5 m during the OA period (black dashed line on Fig. \ref{fig:car_comparison}); and ($ii$) the vehicle does not oscillate during the LK period (in other words, the vehicle moves at a constant speed with $w_f$ and $\theta$ close to $0$ m and $0^\circ$). These feasibility and satisfactory criteria are assumed to be unknown to C-GLISp and are learned based on the expressed feasibility and satisfactory labels.

For the pairwise comparisons, the calibrator expresses her preferences according to the following guidelines: ($i$) whether it is feasible; ($ii$) whether it is satisfactory; ($iii$) whether the vehicle guarantees passengers' comfort during the OA period, for example, by not changing velocities or moving the lateral position too aggressively; ($iv$) whether the deviations of the vehicle velocity from the reference values is minor in both LK and OA periods; ($v$) whether aggressive variations of steering angles are avoided. If a conflict combination among criteria mentioned above appears, criterion ($i$) has the highest priority, and if the conflict is among criteria ($ii$)--($v$), the control policy that leads to qualitatively safer driving practice based on the calibrator's experience is preferred. Note that conditions ($iii$)-($v$) and the method of judging safe driving practice are mainly qualitative/subjective, and it is difficult to express them in terms of quantitative metrics.  

For the example query window illustrated in Fig. \ref{fig:car_comparison}, conflicting combinations of the assessing criteria are observed. Compared to the experiment shown in the left panels of the figure, the experiment shown in the right panels has  more aggressive lateral movements during the OA period. Furthermore, the changes in velocity and steering angle are greater in both frequency and magnitude in both LK and OA periods. The experiment shown in the left panels is feasible since it is implementable ($t_{comp}< T_s$) and stable, while the experiment shown in the right panels of the figure is infeasible since $t_{comp}$ exceeds $T_s$. Additionally, both experiments fail to satisfy the satisfaction conditions. Above all, the performance of the experiment shown in the left panels   is preferred according to criterion ($vi$).

\subsection{Results}
\begin{figure}[t!]
  \centering
  \includegraphics[width=\linewidth]{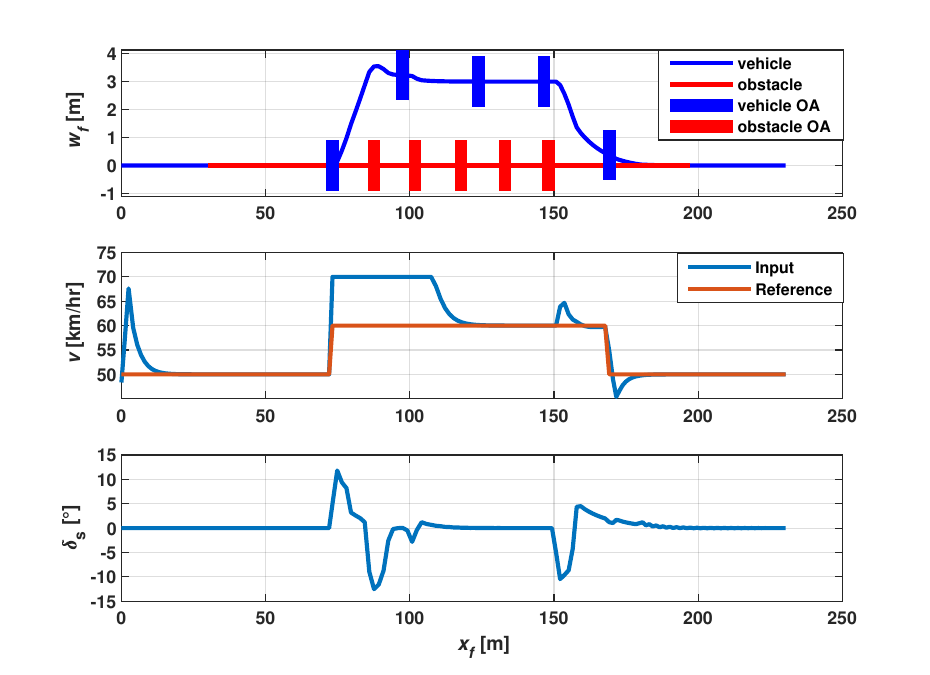}
  \vspace{-0.8cm}
  \caption{Final vehicle control performance obtained by the designed MPC controller. The top subplot shows the vehicle and obstacle positions. The ``vehicle OA'' and ``obstacle OA'' bars show five relative positions of the vehicle and obstacle during the OA period. The middle subplot shows the actual and reference velocity $v$ at different longitudinal positions. The bottom subplot shows the steering angle $\psi$ over the longitudinal position.   }
  \label{fig:car_finalPerf}
  %\vspace{-0.4cm}
\end{figure}

C-GLISp terminates after 50 simulated closed-loop experiments and 49 pairwise comparisons. The best MPC design parameters $T_s$, $\epsilon_c$, $N_p$, $\log(q_{u11})$ and $\log(q_{u22})$ are 0.085 s, 0.100, 23, -0.323 and -3.71, respectively, with a worst-case computation time $t_{comp} =$ 0.0789 s. The closed-loop performance obtained via these MPC design parameters is depicted in Fig.~\ref{fig:car_finalPerf}. As shown in the figure, after only 50 simulated experiments, the proposed algorithm can tune the MPC parameters to achieve feasible and satisfactory performance, accomplishing the driving tasks with smooth and safe maneuvers. 

\section{Conclusions}\label{sec:conclusion_futureExtension}
The algorithm C-GLISp introduced in this paper can handle preference-based
global optimization with unknown objective functions and unknown constraints better than other
existing black-box surrogate methods (PBO and GLISp), as illustrated through benchmark problems. The autonomous driving case study demonstrated the application of C-GLISp in semi-automated MPC calibration. Although convergence to global optimizers cannot be guaranteed, we observed that the C-GLISp can find satisfactory results within a small number of iterations and that it has a higher probability of proposing feasible samples during the exploration thanks to the introduction of additional information by the decision-maker that is used to synthesize corresponding surrogate functions. 

Future research will be devoted to exploring different initial sampling methods, warm-starting procedures in the form of ``transfer learning'' from previous similar optimization runs, and preference-based collective learning to handle multiple decision-makers. We finally note that, while we have used C-GLISp for controller calibration, the algorithm can be used in many other applications in which a few tuning parameters must be decided based on preferences under constraints that cannot be easily quantified.

%\section*{Acknowledgment}

\bibliographystyle{IEEEtran} 
\bibliography{Biblio_pref_unkn_const}    

%\begin{IEEEbiography}
%\end{IEEEbiography}

\end{document}